\documentclass[12pt]{amsart}
\usepackage{latexsym}
\usepackage{graphicx}
\usepackage{amsmath}
\usepackage{amssymb}
\usepackage{verbatim}
\usepackage{amsfonts}

\newtheorem{lemma}{Lemma}
\newtheorem{thm}{Theorem}

\newtheorem{claim}{Claim}

\newcommand{\bd}{\partial}

\newcommand{\h}{\mathbb{H}}

\begin{document}

\title{Curvature bounds for surfaces in hyperbolic 3-manifolds}
\author{William Breslin}
\address{Department of Mathematics, University of Michigan, 530 Church Street, Ann Arbor, MI 48109}
\email{breslin@umich.edu}
\date{February 20, 2007}

\begin{abstract}
A triangulation of a hyperbolic 3-manifold is \textit{L-thick} if each tetrahedron having all vertices in the thick part of $M$ is $L$-bilipschitz diffeomorphic to the standard Euclidean tetrahedron.  We show that there exists a fixed constant $L$ such that every complete hyperbolic 3-manifold has an $L$-thick geodesic triangulation.  We use this to prove existence of universal bounds on the principal curvatures of $\pi_1$-injective surfaces and strongly irreducible Heegaard surfaces in hyperbolic 3-manifolds.
\end{abstract}

\maketitle


\section{Introduction}

Surfaces in a 3-manifold can be complicated for two reasons.  A surface may have small handles and local topology, with regions where the curvature is large.  Also, a surface may be locally well-behaved, but globally fold back on itself in a complicated manner.  It follows from work of Schoen-Yau \cite{Schoen-Yau}, Freedman-Hass-Scott \cite{FHS}, and Schoen \cite{Schoen} that a $\pi_1$-injective embedding (resp. immersion) of a closed surface into a closed hyperbolic 3-manifold can be isotoped (resp. homotoped) to a surface with principal curvatures bounded in absolute value by a universal constant which is independent of the surface and the ambient 3-manifold.  Using different techniques, we extend this to strongly irreducible Heegaard surfaces and to noncompact complete hyperbolic 3-manifolds.

\begin{thm}\label{thm1}
There exists a constant $\omega > 0$ such that the following holds.  If $S$ is a $\pi_1$-injective orientable embedded (resp. immersed) closed surface or a strongly irreducible Heegaard surface in a complete orientable hyperbolic 3-manifold $M$, then $S$ is isotopic (resp. homotopic) to a surface whose principal curvatures are bounded in absolute value by $\omega$.
\end{thm}

These results are new in two ways.  First, the methods of \cite{FHS} do not apply in the presence of geometrically infinite ends or cusps, but Theorem \ref{thm1} holds for any complete hyperbolic 3-manifold.  Second, the curvature bounds of \cite{Schoen} hold only for stable minimal surfaces and Heegaard surfaces are not in general isotopic to stable minimal surfaces.

Theorem \ref{thm1} follows an existence theorem about geodesic triangulations of hyperbolic 3-manifolds.  A geodesic triangulation of a complete hyperbolic 3-manifold $M$ may be forced by the geometry of $M$ to have tetrahedra with very short edges or small dihedral angles.  Big tetrahedra without small dihedral angles cannot live in the thin part of $M$.  In this paper, we show that a complete hyperbolic 3-manifold has a geodesic triangulation such that the tetrahedra contained in the thick part have edge lengths contained in a fixed interval and dihedral angles bounded below by a fixed constant.  The bounds depend only on the constant used to define the thick-thin decomposition of the manifold.  We call such a triangulation a \textit{thick} geodesic triangulation of $M$.

\begin{thm}\label{thm2}
Let $\mu$ be a Margulis constant.  There exist positive constants $a := a(\mu)$, $b := b(\mu)$, and $\theta := \theta (\mu)$ such that the following holds:   A complete hyperbolic 3-manifold $M$ has a geodesic triangulation such that every tetrahedron contained in the thick part of $M$ has dihedral angles bounded below by $\theta$ and edge lengths in $[a,b]$.
\end{thm}

The existence of thick geodesic triangulations holds for higher dimensions as well, which can be found in \cite{Breslin2}.  Emil Saucan has shown that hyperbolic $n$-orbifolds have triangulations whose simplices are uniformly round (called ``fat" triangulations), and he uses this to prove existence of quasi-meromorphic maps which are automorphic with respect to the corresponding Kleinian group (See \cite{saucan1}, \cite{saucan2}, \cite{saucan3}).  However, the triangulations he produces have no uniform bound on the size (i.e. edge lengths) of the simplices, even in the thick part.  The thick triangulations obtained in this paper have uniformly bounded edge lengths, but it should be noted that the lower bound is very small.\\

To prove Theorem \ref{thm2}, we examine Delaunay triangulations of ``well-spaced'' point sets in 3-dimensional hyperbolic space and the problem of eliminating small dihedral angles.  The corresponding question for 3-dimensional Euclidean space has been well-studied.  The only tetrahedra in such a triangulation which can have small dihedral angles are \textit{slivers}, and it was a problem to show how to remove them without creating new ones.  Several techniques for removing slivers have been developed in the Euclidean setting (see \cite{Edels},\cite{Miller},\cite{Li}).  We adapt the technique introduced in \cite{Edels} of perturbing vertices of a Delaunay triangulation in order to remove slivers to the hyperbolic setting.

The paper is organized as follows.  Section 2 contains some definitions and lemmas about geodesic tetrahedra in 3-dimensional hyperbolic space and Section 3 contains the proof of Theorem \ref{thm2}.  We discuss how to use thick triangulations to obtain principal curvature bounds in Section 4.

\textbf{Acknowledgements.}  This work was partially supported by the NSF grant DMS-0135345.  I am grateful to Joel Hass for his guidance and many helpful and stimulating conversations.  I also thank the referee for helpful comments.

\section{Definitions and Lemmas}

\noindent\textbf{Definition.}  Let $t$ be a tetrahedron in $\mathbb{H} ^3$.  Let $p$ be a vertex of $t$.  Let $R_t$ be the circumradius of $t$.  Let $c_p$ be the circumradius of the face opposite $p$.  Let $l_t$ be the length of the shortest edge of $t$.  Let $d_p$ be the distance from $p$ to the plane opposite $p$.  Let $\sigma >0$ , $\rho >0$.  Call $t$ a $(\sigma ,\rho )$\textit{-sliver} if $R_t / l_t \le \rho$ and $d_p / c_p \le \sigma$ for some vertex $p$ of $t$.  Once we have fixed $\sigma$ and $\rho$, we will just refer to $t$ as a \textit{sliver}.\\

The first lemma says that if the edge lengths are bounded between two positive constants and the circumradius is bounded from above, then the only way to have arbitrarily small dihedral angles is to be a sliver.

\begin{lemma}\label{lem1}
If a geodesic tetrahedron $t$ in $\mathbb{H} ^3$ with edge lengths in $[a,b]$ and circumradius at most $R$ has a dihedral angle less than $\theta(a,b,\sigma ) := \arcsin(\frac{\sinh({{\sigma \cdot a / 2}})}{\sinh(b)})$ for some $\sigma >0$, then $t$ is a $(\sigma ,\frac{R}{a} )$-sliver.
\end{lemma}

\begin{figure}\label{fig1}
\includegraphics[width=\textwidth]{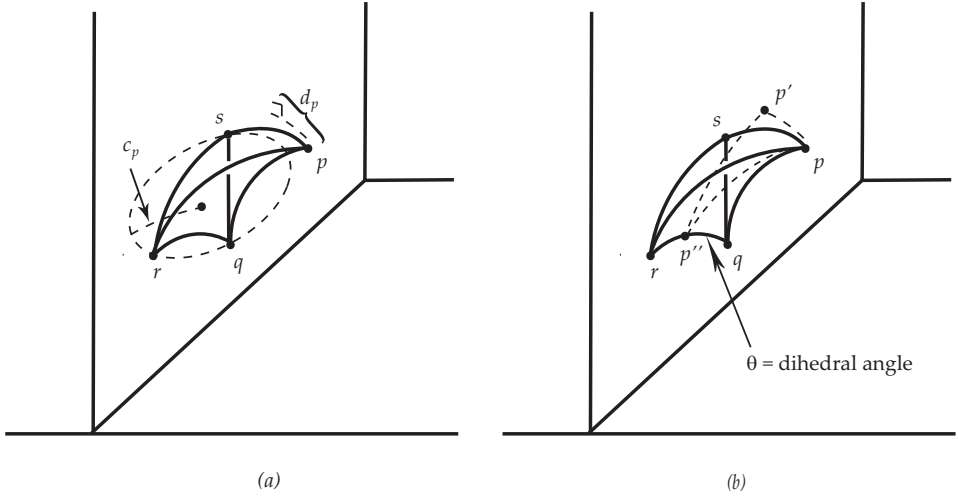}
\caption{(a) The plane containing $[q,r,s]$ is shown as a vertical Euclidean plane in the upper half space model of $\mathbb{H} ^3$.  The circumradius of $[q,r,s]$ is $c_p$.  The distance from $p$ to the plane containing $[q,r,s]$ is $d_p$.  (b) The angle between $[p,r,q]$ and $[q,r,s]$ is $\angle pp''p'$ in the geodesic triangle $[p,p',p'']$.}
\end{figure}

\begin{proof}
Suppose that $t$ has a dihedral angle of less than $\theta$ along the edge $[q,r]$. Project $p$ orthogonally to the plane containing the face $[q,r,s]$ to a point we will call $p'$. Also project $p$ orthogonally to the hyperbolic line containing $[q,r]$, to a point we will call $p''$.  See figure 1.  Now consider the triangle in $\mathbb{H} ^3$ with vertices $p$, $p'$, and $p''$.  The edge between $p$ and $p''$ has length at most $b$.  We have:\\
\begin{center}
${\sinh(||[p,p']||)\over \sin(\theta )}
 = {\sinh(||[p,p'']||)\over \sin({\pi\over 2})}\leq \sinh(b)$\\
\end{center}
Thus
\begin{center}
$||[p,p']|| \leq \operatorname{arcsinh} (\sin(\theta )\cdot \sinh(b))= \sigma\cdot a / 2$
\end{center}
so that
\begin{center}
${{2||[p,p']||}} / a \le \sigma$
\end{center}
The first equality follows from the hyperbolic law of sines \cite{Fenchel}.  Since the lengths of the edges of the triangle $t$ are at least $a$, the circumradius of $[q,r,s]$ is at least $a /2$.  Also, $R_t /l_t \le R / a$.  Thus $t$ is a ($\sigma ,\rho )$-sliver.
\end{proof}

\begin{lemma}\label{lem2}
A geodesic triangle in $\mathbb{H} ^2$ with edge lengths in $[a,b]$ and circumradius at most $R$ has altitudes bounded from below by a constant positive $h_0 := h_0 (a,b,R)$.
\end{lemma}

\begin{figure}\label{fig2}
\includegraphics[width=\textwidth]{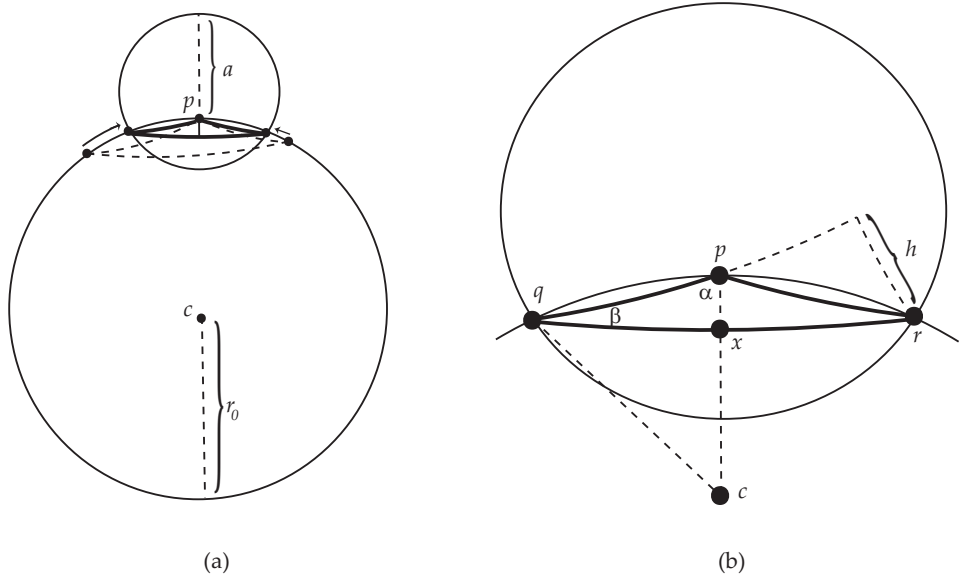}
\caption{(a) If the radius $r_0$ is fixed, then moving $q$ and $r$ closer to $p$ makes the altitude from $p$ smaller.  (b) Both $||[p,x]||$ and $h$ are bounded from below in terms of $a$,$b$, and $R$.}
\end{figure}

\begin{proof}
Since the sum of the angles of $t$ are less than $\pi$ there are at least two angles of $t$ which are less than $\pi / 2$.  Let $p$ be the vertex opposite these angles.  The orthogonal projection of $p$ onto the line containing the opposite edge $[q,r]$ is contained in the interior of $[q,r]$.  Now suppose we have fixed the circumradius $r_0 \in [a / 2 ,R]$ of $t$ and consider all triangles with edges of length at least $a$ such that $p$ projects to the interior of $[q,r]$.  The triangle with the shortest altitude at $p$ is an isosceles triangle which lies on a hyperbolic circle of radius $r_0$ with $||[p,q]|| = ||[p,r]|| = a$ (see figure 2). Let $c$ be the center of the hyperbolic circle containing $p,q,r$.  Let $x$ be the intersection of $[p,c]$ and $[q,r]$.
Let $\beta = \angle pqx$.
Now the altitude of $[p,q,r]$ from $p$ is $||[p,x]||$.  Using the law of cosines, we get
\begin{center}
$\cosh(||[p,x]||) = \cosh(a)\cosh(||[q,x]||) - \sinh(a)\sinh(||[q,x]||)\cos(\beta )$
$\ge \cosh(a)\cosh(||[q,x]||) - \sinh(a)\sinh(||[q,x]||)$.
\end{center}
Let $h_1 (a,r_0 ) = \operatorname{arccosh}(\cosh(a)\cosh(||[q,x]||) - \sinh(a)\sinh(||[q,x]||)$. So far we have shown that the altitude from a vertex which projects to the interior of the opposite face is at least $h_1 (a,r_0 )$ if the circumradius of $[p,q,r]$ is $r_0$.  The triangle for which $\beta$ is minimized and $||[q,x]||$ is maximized has circumradius $R$.  Since $h_1 (a,r_0)$ is minimized when $r_0 = R$, we have that $h_1 (a,R)$ is a lower bound on the altitude from a vertex which projects to the interior of the opposite face for triangles satisfying the hypotheses of the lemma.  Let $h$ be the altitude from $r$.  We have
\begin{center}
$\sin(\beta ) = $ $\frac{\sinh(||[p,x]||)}{\sinh(||[p,q]||)}$ $\ge\frac{\sinh(h_1(a,R))}{\sinh(b)}$.
\end{center}
Also,
\begin{center}
$\sinh(h) =$ $\sinh(||[q,r]||)\sin(\beta ) \ge$ $\sinh(a) \cdot \frac{\sinh(h_1 (a,R))}{\sinh(b)}$,
\end{center}
so that
\begin{center}
$h \ge$ $\operatorname{arcsinh}(\frac{\sinh(a)}{\sinh(b)} \cdot \sinh(h_1 (a,R)))$.
\end{center}
A similar argument works for the altitude from $q$.  Let
\begin{center}
$h_0 (a,b,R) =$ $\operatorname{arcsinh}(\frac{\sinh(a)}{\sinh(b)} \cdot \sinh(h_1 (a,R)))$.
\end{center}
\end{proof}

\begin{lemma}\label{lem3}
If a geodesic tetrahedron $t$ in $\mathbb{H} ^3$ with edge lengths in $[a,b]$ and circumradius at most $R$ is a $(\sigma ,\frac{R}{a} )$-sliver for some $\sigma > 0$, then $d_v /c_v$ is bounded above by a constant $n := n(\sigma ,a,b,R)$ for each vertex $v$ of $t$.  Moreover, $n$ can be chosen so that $n \rightarrow 0$ as $\sigma \rightarrow 0$ and $a,b,R$ remain fixed.
\end{lemma}

\begin{figure}\label{fig3}
\includegraphics[width=\textwidth]{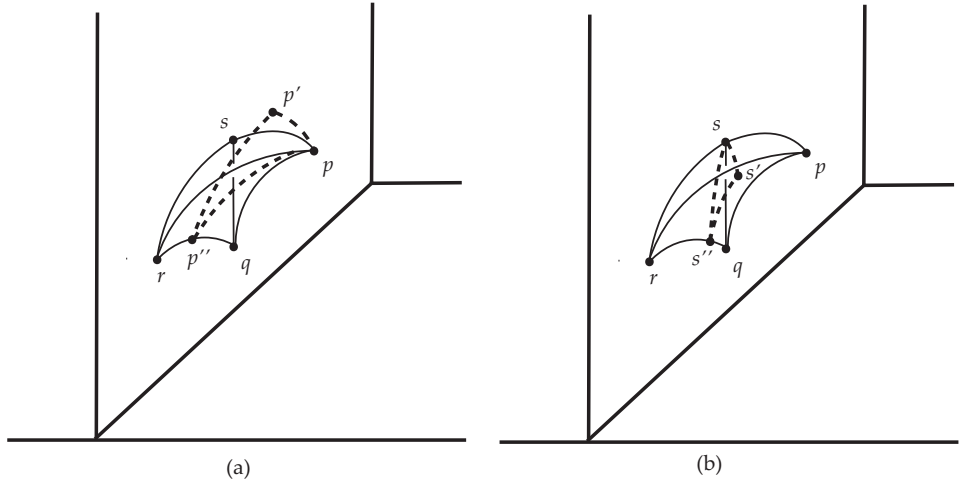}
\caption{(a) To bound $d_s / c_s$, let $p''$ be the projection of $p$ to the line containing $[q,r]$ and consider $[p,p',p'']$.  (b) $\angle ss''s'$ $=$ $\alpha$ $=$ $\angle pp''p'$ is small if $d_p = ||[p,p']||$ is small.}
\end{figure}

\begin{proof}
Since $t$ is a $(\sigma ,\frac{R}{a})$-sliver, we have $d_v /c_v \le \sigma$ for some vertex $v$.  Say $d_p /c_p \le \sigma$.  Let $p'$ be the orthogonal projection of $p$ onto the plane containing $[q,r,s]$.  Let $p''$ be the orthogonal projection of $p$ onto the line containing $[q,s]$.  Let $\alpha = \angle pp''p'$ be the angle between the planes containing $[p,q,r]$ and $[q,r,s]$.    Using the hyperbolic law of sines, we get (see figure 3)
\begin{center}
$\sin(\alpha) = \frac{\sinh(||[p,p']||)}{\sinh(||[p,p'']||)}$.
\end{center}
Now $||[p,p']|| = d_p$ and $||[p,p'']|| \ge h_0 (a,b,R)$, where $h_0 (a,b,R)$ is the constant provided in Lemma \ref{lem2}.  Thus
\begin{center}
$\sin(\alpha) \le \frac{\sinh(d_p )}{\sinh(h_0 (a,b,R))} \le \frac{\sinh(\sigma\cdot R)}{\sinh(h_0 (a,b,R))}$.
\end{center}
Now let $s'$ be the orthogonal projection of $s$ onto the plane containing $[p,q,r]$ and let $s''$ be the orthogonal projection of $s$ onto the line containing $[q,r]$.  Now $\angle ss''s' = \alpha$ and $d_s = ||[s,s']||$, so that
\begin{center}
$\sinh(d_s )) =$ $\sinh(||[s,s'']||) \cdot \sin(\alpha)$  $\le \sinh(b) \cdot  \frac{\sinh(\sigma\cdot R)}{\sinh(h_0 (a,b,R))}$.
\end{center}
Also, the circumradius $c_s$ of $[p,q,r]$ is at least $a\over 2$.  Thus
\begin{center}
$d_s / c_s  \le \frac{2\operatorname{arcsinh}((\sinh(b) \cdot  \frac{\sinh(\sigma\cdot R)}{\sinh(h_0 (a,b,R))})}{a} =: n(\sigma , a,b,R)$.
\end{center}

By projecting $p$ to the other edges of $[q,r,s]$, we can show that $d_q / c_q $ and $d_r / c_r $ are also no more than $n(\sigma , a,b,R)$.
\end{proof}

The next lemma shows that the vertices of a sliver with bounded edge lengths and bounded circumradius all lie near a hyperbolic circle.

\begin{lemma}\label{lem4}
If a geodesic tetrahedron $t = [p,q,r,s]$ in $\mathbb{H} ^3$ with edge lengths in $[a,b]$ and circumradius at most $R$ is a $(\sigma ,\frac{R}{a} )$-sliver for some $\sigma > 0$, then the distance from $p$ to the circumcircle of $[q,r,s]$ is at most an explicit constant $K := K(\sigma ,a,b,R)$.  Moreover, $K$ can be chosen so that $K \rightarrow 0$ as $\sigma \rightarrow 0$ and $a,b,R$ remain fixed.
\end{lemma}

\begin{figure}\label{fig4}
\includegraphics[width=\textwidth]{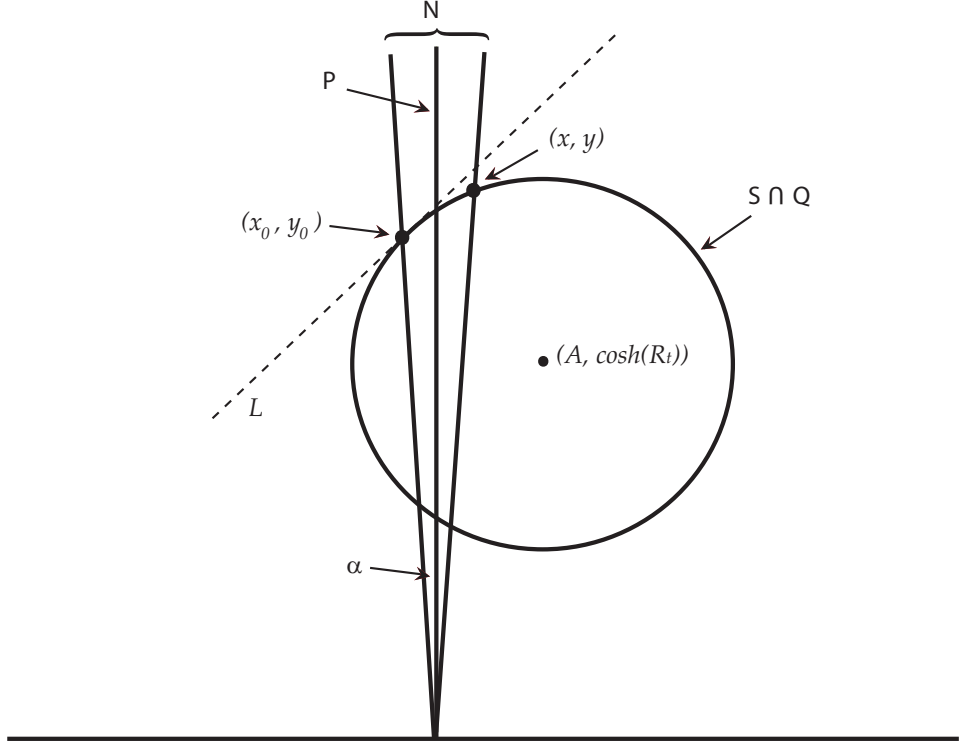}
\caption{This shows the circle $S\cap\mathcal{Q}$ in the plane $\mathcal{Q}$. The plane $\mathcal{P}$ contains $[q,r,s]$ and $p$ is either $(x_0 ,y_0 )$ or $(x, y)$.  $\mathcal{N}$ is the $d(p, \mathcal{P})$-neighborhood of $\mathcal{P}$.  $L$ is the Euclidean line in $\mathcal{Q}$ through $(x_0 ,y_0)$ which is tangent to $S\cap\mathcal{Q}$.}
\end{figure}

\begin{proof}
Let $S$ be the circumsphere of $[p,q,r,s]$.  Let $R_t$ be the radius of $S$.  Let ${C}$ be the circumcircle of $[q,r,s]$. Let $c$ be the radius of ${C}$.  We will use an upper-half space model for $\mathbb{H} ^{3}$ such that ${C}$ is contained in the plane $\mathcal{P}=\{ x=0\}$, the smaller component of $S - \mathcal{P}$ is contained in $\{ x\leq 0\}$, the hyperbolic center of $S$ has coordinates $(A,0,1)$, $p$ is contained in the plane $\mathcal{Q} = \{ y=0\}$, and the $z$-coordinate of $p$ is $\ge 1$.  Figure 4 shows the plane $\mathcal{Q}$.  From now on, we will use $(x,z)$-coordinates when working in $\mathcal{Q}$.  Since $S\cap\mathcal{Q}$ is a hyperbolic circle of radius $R_t$ centered at $(A,1)$, it is a Euclidean circle of radius $\sinh (R_t)$ centered at $(A,\cosh (R_t))$.  We can express $A$ as a function of $R_t$ and $c$:

\begin{center}
$A=A(R_t,c)=\sqrt{(\sinh({R_t}))^2 - (\sinh(c))^2  }$,
\end{center}

\noindent and describe $S\cap\mathcal{Q}$ by the equation

\begin{center}
$(x-A)^{2} + (y-\cosh(R_t))^{2} = (\sinh(R_t))^{2}$.
\end{center}

Let $w$ be the orthogonal projection of $p$ to $\mathcal{P}$.  Since $[p,q,r,s]$ is a sliver, we have
\begin{center}
$d(p,w)\le n(\sigma ,a,b,R)\cdot c\leq n(\sigma ,a,b,R) \cdot R$,
\end{center}
where $n(\sigma ,a,b,R)$ is the constant from Lemma \ref{lem3}.

Now we will find an upper bound for $d(w,{C}) = d(w,u)$, where $u = S\cap\mathcal{Q}\cap\mathcal{P}$.  Let $\beta = \angle wup$.  We have
\begin{center}
$\sin(\beta) = \frac{\sinh(c)}{\sinh(R_t)} \ge \frac{\sinh(a/2)}{\sinh(R)}$.
\end{center}
Let $\mathcal{N}$ be the $d(p,\mathcal{P})$-neighborhood of $x=0$ in $\mathcal{Q}$.  Now $\mathcal{N}$ is the region in $y\ge 0$ between the Euclidean lines $y= x/{\tan(\alpha)}$ and $y= -x/{\tan(\alpha)}$, where $\alpha = {\pi - 4{\arctan} (e^{-d(p,\mathcal{P})})}$ is the angle between each line and the $y$-axis.  Note that $\alpha \le J(\sigma,a,b,R) := [\pi - 4{\arctan} (e^{-n(\sigma ,a,b,R) \cdot R})]$ and $J(\sigma,a,b,R) \rightarrow 0$ as $\sigma \rightarrow 0$.  We have
\begin{center}
$d_E (w,p) = \tan(\alpha)\cdot d_E (w,\{ z=0\}) \le \tan(J(\sigma,a,b,R)) \cdot 2\cosh(R)$,
\end{center}
where $d_E$ denotes Euclidean distance.
We also have
\begin{align*}
d_E (w,u) &= \frac{d_E (w,p) \cos(\beta)}{\sin(\beta)} \\
&\le \frac{2\tan(J(\sigma,a,b,R))\cosh(R)\sinh(R)}{\sinh(a/2)}.
\end{align*}
Since the $z$-coordinates of $w$ and $u$ are at least $1$, the upper bound on $d_E (w,u)$ gives us an upper bound on $d(w,u)$:
\begin{center}
$d(w,u) \le \frac{2\tan(J(\sigma,a,b,R))\cosh(R)\sinh(R)}{\sinh(a/2)}$.
\end{center}
We now have
\begin{align*}
d(p,C) &\le d(p,w) + d(w,u)\\
&\le n(\sigma ,a,b,R) \cdot R + \frac{2\tan(J(\sigma,a,b,R))\cosh(R)\sinh(R)}{\sinh(a/2)}.
\end{align*}
Let $K(\sigma,a,b,R) := n(\sigma ,a,b,R) \cdot R + \frac{2\tan(J(\sigma,a,b,R))\cosh(R)\sinh(R)}{\sinh(a/2)}$.
\end{proof}

\noindent\textbf{Definition.}  Let $T =[q,r,s]$ be a geodesic triangle in $\mathbb{H} ^3$ with edge lengths in $[a,b]$ and circumradius at most $R$.  Let $\sigma > 0$.  The $(\sigma ,a,b,R)$-\textit{sliver region of} $T$ is the set of points $p$ in $\mathbb{H} ^3$ such that $[p,q,r,s]$ is a $(\sigma ,\frac{R}{a} )$-sliver with edge lengths in $[a,b]$ and circumradius at most $R$.  We will denote this region by $sliver_{(\sigma ,a,b,R)}[q,r,s]$.

\begin{lemma}\label{lem5}
The volume of the $(\sigma ,a,b,R)$-sliver region of a geodesic triangle $T = [q,r,s]$ in $\mathbb{H} ^3$ with edge lengths in $[a,b]$ and circumradius at most $R$ is at most an explicit constant $V := V(\sigma ,a,b,R)$.  Moreover, $V$ can be chosen so that $V \rightarrow 0$ as $\sigma \rightarrow 0$ and $a,b,R$ remain fixed.
\end{lemma}

\begin{figure}\label{fig5}
\includegraphics[width=\textwidth]{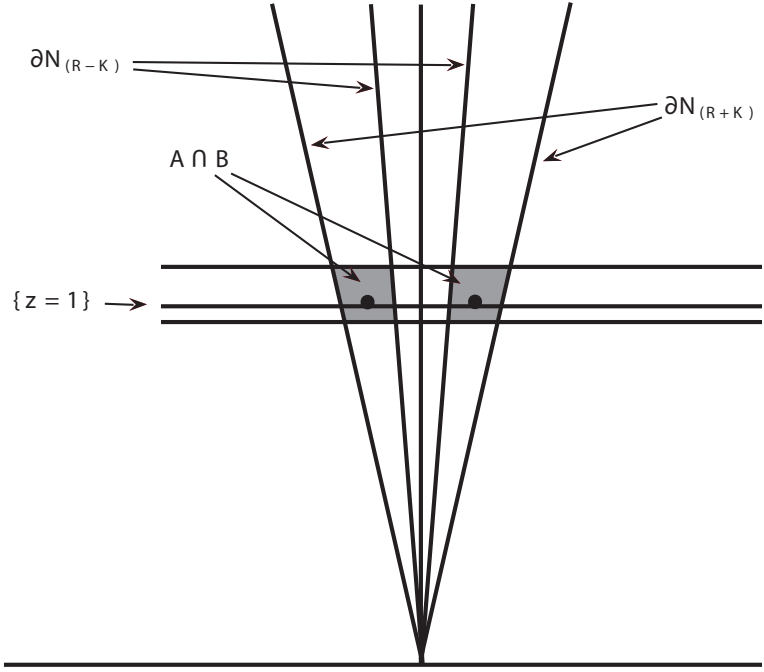}
\caption{The dots are the circumcircle of $[q,r,s]$, which has radius $R$.  The shaded region is $A\cap B$, which contains the vertex $p$.}
\end{figure}

\begin{proof}
It follows from Lemma \ref{lem4} that $sliver_{(\sigma ,a,b,R)}[q,r,s]$ is contained in a $K$-neighborhood, $U$, of the circumcircle of $[q,r,s]$, where $K := K(\sigma ,a,b,R)$ is the constant provided by Lemma \ref{lem4}. Let $V := V(\sigma ,a,b,R)$ be the volume of $U$, so that $V$ is an upper bound on the volume of $sliver_{(\sigma ,a,b,R)}$.  Since $K \rightarrow 0$ as $\sigma \rightarrow 0$, we have that $V \rightarrow 0$ as $\sigma \rightarrow 0$.
\end{proof}

\section{Thick geodesic triangulations}

We will prove Theorem \ref{thm2} in this section.  Let $M$ be a complete orientable hyperbolic 3-manifold.\\

\noindent\textbf{Definition.}  Let $\mu >0$.    The $\mu$-\textit{thick part} of $M$, denoted by $M_{[\mu ,\infty )}$ is the set of points where the injectivity radius is at least $\mu / 2$.  The $\mu$-\textit{thin part} of $M$, denoted by $M_{(0, \mu ]}$, is the closure of the complement of $M_{[\mu, \infty)}$.\\

The Margulis lemma \cite{Margulis} implies that there exists a constant $\mu_0$ (a \textit{Margulis constant}) such that for $\mu\le\mu_0$ every component of the $\mu$-thin part of any complete orientable hyperbolic 3-manifold is either a tubular neighborhood of a closed geodesic or the quotient of a horoball by an abelian parabolic subgroup.  We will refer to a component of the $\mu$-thin part which is a tubular neighborhood of a closed geodesic as a $\mu$\textit{-tube}.  Let $\mu\le\mu_0$.

Let $\mathcal{S}$ be a generic set of points in $M$ such that for any $p$ $\in$ $M$ the ball $B(p,$inj$(M,p)/5)$ centered at $p$ with radius inj$(M,p)/5$ contains a point of $\mathcal{S}$ in its interior. We require the points to be denser where the injectivity radius is smaller so that the Delaunay triangulation of $\mathcal{S}$ is defined.\\

\noindent\textbf{Definition.} The \textit{Delaunay triangulation of $\mathcal{S}$} is the geodesic triangulation of $M$ determined as follows:  A set, $\{ p,q,r,s \}$, of four vertices in $\mathcal{S}$ determines a tetrahedron in $\mathcal{T}$ if and only if the minimal radius circumscribing sphere contains no points of $\mathcal{S}$ in its interior.\\

See \cite{Leibon} for existence of Delaunay triangulations in Riemannian manifolds.
One disadvantage of using a Delaunay triangulation is that without further assumptions on the set of points, the tetrahedra may have arbitrarily small dihedral angles.  We now add a further restriction on the set $\mathcal{S}$ which will restrict the way a tetrahedron which intersects the thick part can have small dihedral angles.  Assume that $\mathcal{S}$ is maximal with respect to the condition that for each point $p$ in $\mathcal{S} \cap M_{[\mu ,\infty)}$ we have $d(p,q) \ge \epsilon := \mu /100$ for every $q$ in $\mathcal{S} \setminus \{ p \}$.  We have defined $\epsilon$ to be sufficiently small with respect to the injectivity radius so that if we perturb a point $p\in\mathcal{S}$, any changes to the Delaunay triangulation of $\mathcal{S}$ will occur in a ball which lifts to the universal cover.  Thus we may work in $\mathbb{H} ^3$.

This extra assumption on $\mathcal{S}$ also implies that a tetrahedron $t$ in the Delaunay triangulation of $\mathcal{S}$ has circumradius at most $\epsilon$, and edge lengths in the interval $[\epsilon ,2\epsilon ]$.  Thus such a tetrahedron has small dihedral angles only if it is a sliver (Lemma \ref{lem1}) and such a tetrahedron is a sliver only if its vertices are close to a circle (Lemma \ref{lem4}).\\

Let $\delta := \delta(\mu) = \epsilon (\mu) /10$.  \\

\noindent\textbf{Definition.} A \textit{good perturbation} of $\mathcal{S}$ is a collection of points $\mathcal{S'}$ in $M$ such that there exists a bijection $\phi : \mathcal{S} \rightarrow \mathcal{S'}$ with $d(p,\phi(p)) \le \delta$ for every $p \in \mathcal{S}$. Denote $\phi(p)$ by $p'$.  If $\mathcal{T}$ and $\mathcal{T'}$ are the Delaunay triangulations of $\mathcal{S}$ and $\mathcal{S'}$, then we will say that $\mathcal{T'}$ is a \textit{good perturbation} of $\mathcal{T}$.\\

Note that if $t$ is a tetrahedron in the Delaunay triangulation of a good perturbation of $\mathcal{S}$ which is contained in $M_{[\mu ,\infty )}$, then $t$ has circumradius no more than $\epsilon + \delta$ and edge lengths between $\epsilon - 2\delta$ and $2\epsilon + 2\delta$.  We can still apply our lemmas to such tetrahedra.

\begin{lemma}\label{lem6}
Let $p\in\mathcal{S} \cap M_{[\mu ,\infty)}$.  The number of triples, $\{ q,r,s\}$, of points in $\mathcal{S}$ such that $[p',q',r',s']$ is a tetrahedron in some good perturbation of $\mathcal{T}$ is bounded from above by a constant $N := N(\mu )$.
\end{lemma}

\begin{proof}
The $(\frac{\epsilon}{2} -\delta)$-balls centered at the points of $\mathcal{S} \cap M_{[\mu ,\infty)}$ are mutually disjoint since no two points of $\mathcal{S}  \cap M_{[\mu ,\infty)}$ are closer than $\epsilon - 2\delta$ to each other.  If $p'$ and $q'$ are vertices of a tetrahedron in the Delaunay triangulation $\mathcal{T'}$ of $\mathcal{S'}$, then $d(p',q') \le 2\epsilon + 2\delta$, so that the $(\frac{\epsilon}{2} -\delta)$-ball centered at $q'$ is contained in the $(2\epsilon + 2\delta)$-ball centered at $p'$.  There can be at most
\begin{center}
$m := m(\mu ) = [\frac{vol(B(2\epsilon + 2\delta))}{vol(B(\frac{\epsilon}{2} -\delta))}]$
\end{center}
mutually disjoint $(\frac{\epsilon}{2} -\delta)$-balls contained in a $(2\epsilon + 2\delta)$-ball, where $[w]$ is the integer part of $w$.  One of these is centered at $p'$.  So there are at most $m - 1$ vertices in $\mathcal{S}$ which may be the vertex of a tetrahedron in $\mathcal{T'}$ which also has $p'$ as a vertex.  Thus the number of triples $\{ q,r,s \}$ of points in $\mathcal{S}$ such that $[p',q',r',s']$ is a tetrahedron in some good perturbation of $\mathcal{T}$ is at most $m\choose{3}$.  Let $N(\mu ) := {m\choose{3}}$.
\end{proof}

\textit{Proof of Theorem \ref{thm2}.}  Our proof is based on a method introduced in \cite{Edels} to remove slivers from triangulations of $\mathbb{E} ^3$.  The plan is to perturb one point of $\mathcal{S}\cap\ M_{[\mu,\infty)}$ at a time as follows.  Let $\rho = (\epsilon + \delta ) / (\epsilon - 2\delta )$.   Let $\sigma$ be a positive constant to be determined later.  Let ${p_{1}}\in\mathcal{S}\cap\ M_{[\mu,\infty)}$.  Let ${\mathcal{U}}_{1}$ be the set triangles $[q,r,s]\in\mathcal{T}$ such that there exists a good perturbation $\mathcal{T'}$ of $\mathcal{T}$ which is obtained by perturbing only the point $p_{1}$ and such that $[{p_{1}}',q,r,s]\in\mathcal{T'}$.  We want to pick ${p_{1}}'$ in the ball of radius $\delta$ centered at $p_{1}$ (so that the perturbation is good) and outside the $(\sigma ,\rho )$-sliver region of every triangle in ${\mathcal{U}}_{1}$.  Assume that we can find such a point ${p_1}'$ and call the new set of points ${\mathcal{S}}_{1}$ and the new triangulation ${\mathcal{T}}_{1}$.  Assume we have perturbed the points $p_{1},...,p_{n}$ to ${p_{1}}',...,{p_{n}}'$ and now have a set of points ${\mathcal{S}}_{n}$ and a triangulation ${\mathcal{T}}_{n}$ such that none of ${p_{1}}',...,{p_{n}}'$ is the vertex of a sliver. Let $p_{n+1}$ be a point in $[{\mathcal{S}}_{n} \cap M_{[\mu,\infty)}]     -\{{p_{1}}',...,{p_{n}}'\}$.  Let ${\mathcal{U}}_{n+1}$ be the set of triangles $[q,r,s]\in{\mathcal{T}}_{n}$ such that there exists a good perturbation ${\mathcal{T}_{n}}'$ of ${\mathcal{T}}_{n}$ which is obtained by perturbing only the point $p_{n+1}$ and such that $[{p_{n+1}}',q,r,s]\in{{\mathcal{T}_{n}}'}$.  We choose a point ${p_{n+1}}'$ in the ball of radius $\delta$ centered at $p_{n+1}$ and outside the $(\sigma ,\rho )$-sliver region of every triangle in ${\mathcal{U}}_{n+1}$.

Suppose $M$ has finite volume.  Let $\mathcal{T'}$ be the triangulation we get after perturbing every point of $\mathcal{S} \cap M_{[\mu ,\infty )}$ once and only once (There are only finitely many since $M$ has finite volume).  Let $[p',q',r',s']\in\mathcal{T'}$.  Suppose that $p'$ was the last point perturbed among these four points.  We chose $p'$ to be outside $sliver[q',r',s']$, so that $[p',q',r',s']$ is not a sliver.  Thus any tetrahedron of $\mathcal{T'}$ contained in $M_{[\mu, \infty)}$ is not a sliver.

If $M$ has infinite volume, then the above procedure can be used to perturb the vertices contained in an $N$-ball centered at some fixed point $x_0$, giving us a geodesic triangulation $\mathcal{T}_N$ of $M$ such that any tetrahedron contained in $M_{[\mu, \infty)} \cap B(x_0,N)$ is not a sliver.  Suppose we want to define the final triangulation on $B(x_0,N)$.  Since the triangulations $\mathcal{T}_n$ can be chosen to agree on the ball $B(x_0,N)$ for $n\ge 100N$, we can use the triangulation $\mathcal{T}_{100N}$ to define the triangulation inside $B(x_0,N)$.

We will now prove that for suitable $\sigma$, it is always possible to find a point within $\delta$ of the original point and outside the concerned $(\sigma ,\rho )$-sliver regions.

Suppose we are considering the point $p$ at some step in our process.  We will use $a = \epsilon - 2\delta$, $b = 2\epsilon + 2\delta$, and $R = \epsilon + \delta$ in what follows.  Any tetrahedron in a good perturbation of $\mathcal{T}$ which is contained in $M_{[\mu ,\infty )}$ has edge lengths in $[a,b]$ and circumradius no more than $R$.  The total volume of all the sliver regions we need to consider when perturbing $p$ is at most $N(\mu )*V(\sigma ,a,b,R)$, where $N(\mu )$ is the constant from Lemma \ref{lem6} and $V(\sigma ,a,b,R)$ is the constant from Lemma \ref{lem5}.  Since $V(\sigma ,a,b,R) \rightarrow 0$ as $\sigma \rightarrow 0$, we may choose $\sigma > 0$ so small that the total volume of the sliver regions under consideration is less than the volume of the $\delta$-ball centered at $p$.  This implies that we may find a point $p'$ in the $\delta$-ball centered at $p$ which is not contained in any of the concerned sliver regions.  \hfill $\Box$\\

Note that the tetrahedra in a triangulation provided by Theorem \ref{thm2} which are contained in the thick part of the manifold will be $L$-bilipschitz diffeomorphic to the standard Euclidean tetrahedron for a fixed constant $L$.

\section{Obtaining principal curvature bounds}

We begin this section with some definitions and results from normal surface theory.\\

\noindent\textbf{Definition.} A \textit{normal arc} on the face of a tetrahedron is a properly embedded arc with endpoints in distinct edges of the face.  A \textit{normal curve} on the boundary of a tetrahedron is an embedded closed curve which is transverse to the 1-skeleton and whose intersection with each face is a collection of normal arcs.  A \textit{normal disk} in a tetrahedron is a properly embedded disk whose boundary is a normal curve which intersects the 1-skeleton in 3 or 4 points.  An \textit{almost normal piece} is either a properly embedded disk whose boundary is a normal curve which meets the 1-skeleton in 8 points or two normal disks joined by an unknotted tube. \\

In a given tetrahedron, there are finitely many normal disks and almost normal pieces, up to isotopy preserving the faces, edges, and vertices.\\

\noindent\textbf{Definition.} Let $M$ be an orientable 3-manifold with a triangulation $\mathcal{T}$.  A closed surface $S$ embedded in $M$ is \textit{normal} with respect to $\mathcal{T}$ if the intersection of $S$ with each tetrahedron in $\mathcal{T}$ in a union of normal disks.  We say an embedded surface $S$ is \textit{almost normal} with respect to $\mathcal{T}$ if the intersection of $S$ with one tetrahedron of $\mathcal{T}$ is a union of one almost normal piece and a collection of normal disks and the intersection of $S$ with every other tetrahedron is a union of normal disks. An immersion $f: S \rightarrow M$ is \textit{normal} with respect to $\mathcal{T}$ if for each tetrahedron $t$ in $\mathcal{T}$ every connected component $D$ of $f^{-1}(t)$ is a disk on which $f$ is injective and $f(D)$ is a normal disk.\\

\noindent\textbf{Theorem A} (\cite{Haken}, \cite{Casson})  A closed $\pi_1$-injective surface embedded (resp. immersed) in a triangulated irreducible 3-manifold can be isotoped (resp. homotoped) to a surface which is normal with respect to the given triangulation.\\

\noindent\textbf{Theorem B} (\cite{Rub} , \cite{Stocking}) A strongly irreducible Heegaard surface in a triangulated irreducible 3-manifold can be isotoped to a surface which is almost normal with respect to the given triangulation.\\

In this section, we prove existence of universal bounds on the principal curvatures of $\pi_1$-injective immersed surfaces and strongly irreducible Heegaard surfaces in hyperbolic 3-manifolds (i.e Theorem \ref{thm1}).  We will use a thick geodesic triangulation provided by Theorem \ref{thm2} and normal surface theory.   Normal surface theory implies that the surfaces can be put into normal or almost normal form with respect to any topological triangulation.  To obtain bounds on principal curvatures we assume that the triangulation is thick. Then we choose specific normal and almost normal pieces for each tetrahedron in $\h^3$ so that the pieces meet to form a smooth surface when two tetrahedra are glued together via an isometry, and so that the principal curvatures of the normal and almost normal pieces vary continuously as the vertices of the tetrahedron vary. This will provide a universal bound on the principal curvatures of normal or almost normal surfaces, since tetrahedra from a thick triangulation in the thick part of a hyperbolic 3-manifold come from a fixed compact set.\\

\noindent\textbf{\textit{Geometric pieces.}}
Theorem \ref{thm2} implies that each hyperbolic 3-manifold $M$ has a geodesic triangulation such that each tetrahedron contained in the thick part of $M$ has dihedral angles and edge lengths bounded by fixed constants.  The tetrahedra of such a triangulation contained in the thick part of $M$ are all $L$-bilipschitz diffeomorphic to the standard Euclidean tetrahedron for a fixed constant $L$.  Thus each hyperbolic 3-manifold has a geodesic triangulation such that the tetrahedra contained in the thick part of $M$ are contained in a fixed compact set of tetrahedra which depends only on the constant $L$.
For a positive constant $L$, let $\mathcal{R} (L)$ be the set of geodesic tetrahedra in $\h^3$ which are $L$-bilipschitz diffeomorphic to the standard Euclidean tetrahedron considered up to isometry. Fix a parametrization of $\mathcal{R} (L)$.    For each geodesic triangle which appears as the face of a tetrahedron in $\mathcal{R} (L)$, choose normal arcs of each type which are smooth, meet each edge orthogonally at the midpoint of the edge, and meet each other only at the midpoints of the edges.  Choose these normal arcs so that they change continuously as the vertices of the triangle are moved.   Choose $\rho > 0$ so small that for each chosen arc, we can choose two more smooth arcs of the same type (one on each side) which meet the edges orthogonally a distance $3\rho /2$ from the midpoint and so that the distance of each point of a new arc to the original arc is between $\rho$ and $2\rho$.  Again, choose these arcs so that they change continuously as the vertices of a triangle are moved.

Now we will choose normal and almost normal pieces of each type which are bounded by the normal arcs chosen above.  Let $t\in\mathcal{R} (L)$.  First choose normal pieces of each type in $t$ and an almost normal octagon which are smooth, bounded by the original normal arcs chosen above (for the corresponding faces), meet the boundary of $t$ orthogonally, and meet each other only at their boundaries (except for pairs of quads and octagons).  Choose these normal disks so that they change continuously as the vertices of the tetrahedron are moved.  Note that we have chosen three pieces for each type of normal or almost normal piece. Also tube pairs of normal disks together to form the other almost normal pieces, so that they also change continuously as the shape of a tetrahedron changes.  We will call these normal and almost normal pieces $\textit{geometric pieces}$.  See figure 6.  We may assume that the geometric pieces are ruled annuli in a small neighborhood of each face.  Since $\mathcal{R} (L)$ is compact and there are finitely many geometric pieces in each tetrahedron, the principal curvatures of the geometric pieces are uniformly bounded in absolute value.  Let $B(L)$ be an upper bound on the absolute value of the principal curvature of the geometric pieces for tetrahedra in $\mathcal{R} (L)$.  Given an embedded normal surface with respect to a triangulation with tetrahedra from $\mathcal{R} (L)$, we can isotop each normal and almost normal piece (keeping it embedded) to be so close (in the $C^{\infty}$ topology) to one of the corresponding geometric pieces that the principal curvatures are between $-(B(L) + 1)$ and $B(L) + 1$.\\

\begin{figure}\label{fig6}
\includegraphics[width=\textwidth]{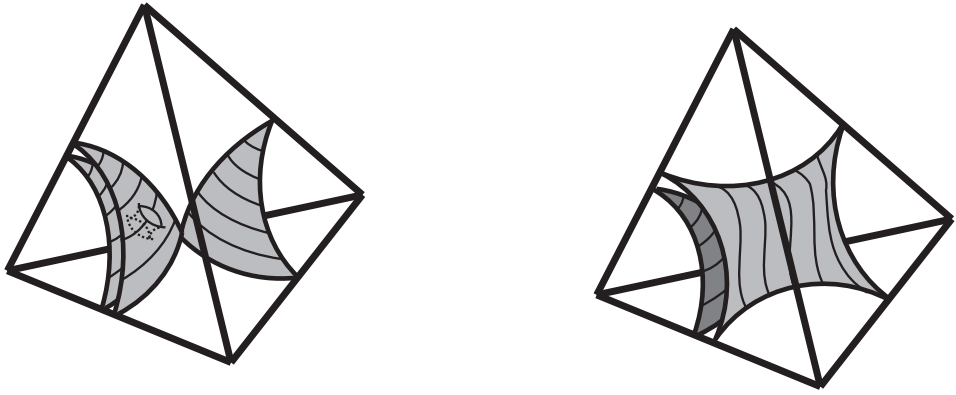}
\caption{Some geometric normal and almost normal pieces.}
\end{figure}

Using geometric pieces in a thick triangulation and normal surface theory (i.e. Theorem A and Theorem B), we can bound the principal curvatures of  $\pi_1$-injective surfaces and strongly irreducible Heegaard surfaces in the thick part of a hyperbolic 3-manifold.  The final step is to extend these bounds to the thin part.  This is easy for Heegaard surfaces because the surface can be isotoped into the thick part and normalized in the thick part.  More work is required in the case of a $\pi_1$-injective surface.\\

\textit{Proof of Theorem \ref{thm1} when $S$ is a strongly irreducible Heegaard surface.} Let $S$ be a strongly irreducible Heegaard splitting of a complete orientable hyperbolic 3-manifold $M$.
Using Theorem A from the introduction we can choose $0 < \epsilon _0 < \epsilon _1 << 1$ so that the following holds:  First, $\epsilon_1$ is less than the 3-dimensional Margulis constant.  Also, if $T_0$ is an $\epsilon _0$-tube in a hyperbolic 3-manifold and $T_1$ is the corresponding $\epsilon _1$-tube, then the distance from $\partial T_0$ to $\partial T_1$ is at least $1$. Let $\mathcal{T}$ be a geodesic triangulation of $M$ provided by Theorem \ref{thm2} with $\mu = \epsilon_0$.
Isotope $S$ so that it is contained in $M_{[\epsilon _1 ,\infty)}$.  We can do this by first isotoping a spine for one of the handlebodies into the interior of $M_{[\mu, \infty )}$ and then isotoping $S$ to the boundary of a small neighborhood of the spine.
The length of an edge of a tetrahedron which has a vertex in $M_{(0 ,\epsilon _0 ]}$ is less than $1$ (since $\epsilon _1 << 1$), so that $S$ intersects only tetrahedra with vertices in $M_{[\epsilon _0 ,\infty )}$.  Theorem B implies that $S$ is isotopic to surface which is almost normal with respect to $\mathcal{T}$.  Also, the normalizing procedure used in \cite{Stocking} does not push $S$ into any tetrahedra which $S$ does not already intersect.  Thus $S$ is isotopic to an almost normal surface which intersects only tetrahedra having all vertices in $M_{[\epsilon _0,\infty)}$.  Each such tetrahedron is isometric to a tetrahedra in $\mathcal{R}(L)$.  Isotope each normal and almost normal piece of $S$ to be so close to the corresponding special normal disk that it has principal curvatures bounded in absolute value by $B(L) + 1$. \hfill $\Box$\\

\textit{Proof of Theorem \ref{thm1} when $S$ is $\pi_1$-injective.} Choose $\epsilon _0 < \epsilon _1 < \epsilon _2 < \epsilon _3 << 1$ so that the following holds.  First, $\epsilon _3$ is less than the Margulis constant.  Also, if $T_0$ is an  $\epsilon _0$-tube in a complete hyperbolic 3-manifold and $T_1$, $T_2$, $T_3$ are the corresponding $\epsilon _1$-, $\epsilon _2$-, and $\epsilon _3$- tubes, then $d(\partial T_0 ,\partial T_1 )$, $d(\partial T_1 , \partial T_2 )$ and $d(\partial T_2 , \partial T_3 )$ are all at least 1.  Let $\mathcal{T}$ be a geodesic triangulation of $M$ provided by Theorem \ref{thm2} with $\mu = \epsilon _0$.

Assume $f: S \rightarrow M$ is an embedding.  Since $S$ is incompressible and closed, we can isotop it so that it misses any $\epsilon _3$-cusps and meets each $\epsilon _3$-tube in either the empty set or a finite set of disjoint totally geodesic disks orthogonal to the boundary of the tube.  Assume that $S$ meets each $\epsilon _3$-tube in as few disks as possible.  By Theorem A we can normalize $S$ with respect to $\mathcal{T}$.

\begin{claim}\label{claim1}
After normalization, $S$ still meets each $\epsilon_2$-tube in either the empty set or a finite set of disjoint totally geodesic disks orthogonal to the boundary of the tube.
\end{claim}

\begin{figure}\label{fig7}
\includegraphics[width=\textwidth]{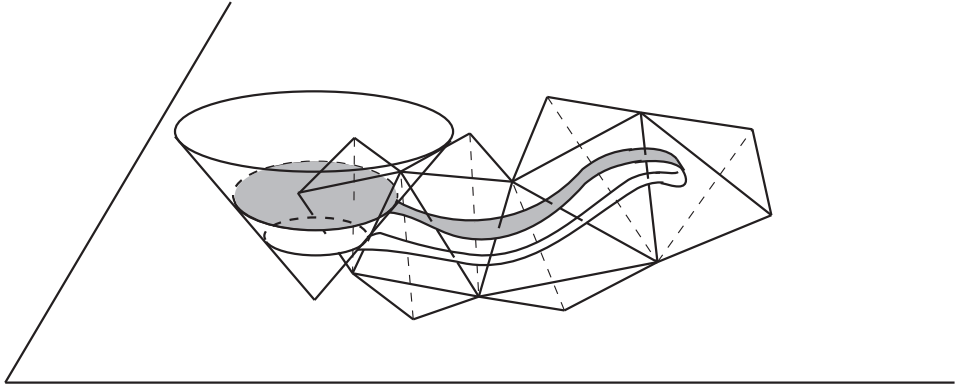}
\caption{A lift of the tube $T_2$, part of $S$, and some tetrahedra to $\h^3$.}
\end{figure}

\begin{proof}
Since the distance between $\partial M_{(0,{\epsilon _2})}$ and $\partial M _{(0,{\epsilon _3})}$ is at least $1$ and the edge of any tetrahedron in $\mathcal{T}$ which meets $M_{[{\epsilon _3}, \infty)}$ is less than $1$, any tetrahedron in $\mathcal{T}$ which meets $M_{[{\epsilon _3}, \infty )}$ is disjoint from $M_{[{\epsilon_2},\infty)}$.  Also, normalization never pushes a surface into a tetrahedron previously disjoint from the surface.  Thus our normalized surface does not meet any ${\epsilon}_2$-cusps or ${\epsilon}_2$-tubes which it did not meet before.
Now consider one of our totally geodesic disks $D$, which meets an ${\epsilon}_3$-tube $T_3$.  The intersection of $D$ with a totally geodesic tetrahedron $t$ is a normal disk if $D\cap t$ is contained in the interior of $D$.  So $D$ intersects any tetrahedron $t$ which meets the corresponding ${\epsilon}_2$-tube,$T_2$, in a normal disk.  Suppose that we are forced to isotop $D\cap T_2$ after some finite number of steps.  A totally geodesic disk cannot intersect the interior of a face in a trivial curve, so the isotopy which moves $D\cap T_2$ is an isotopy which eliminates an arc in $S\cap f$ which intersects an edge of $f$ twice, for some face $f$ of a tetrahedron which meets $T_2$.  See figure 7.  This implies that $D$ is connected to another totally geodesic disk in $T_2 \cap S$ by a strip which is isotopic into $\partial {T_2}$.  This contradicts the assumption that $S$ meets each ${\epsilon}_3$-tube in as few disks as possible.  Thus $S$ still meets each $\epsilon_2$-tube in either the empty set or a set of totally geodesic disks after normalization.
\end{proof}

Let $S'$ be the surface obtained from $S$ by isotoping the normal pieces of $S$ to be so close to the corresponding geometric piece that the curvature of $S'$ is bounded by $B(L) + 1$ in the $\epsilon_2$-thick part of $M$.

We now have two surfaces $S$ and $S'$.  The surface $S$ is totally geodesic in the $\epsilon_2$-thin part so its principal curvatures are bounded in the $\epsilon_2$-thin part, but we only know that $S$ is a normal surface in the $\epsilon_2$-thick part.  The principal curvatures of $S'$ are bounded by $B(L) + 1$ in the $\epsilon_2$-thick part, but the geometric pieces in the thin part must have points where the principal curvatures are large.  The last step is to show that $S\cap M_{(0,\epsilon_2]}$ can be joined to $S'\cap M_{[\epsilon_2, \infty)}$ to get a surface isotopic to $S$ with universally bounded principal curvatures.  See figure 8.

\begin{figure}\label{fig8}
\includegraphics[width=\textwidth]{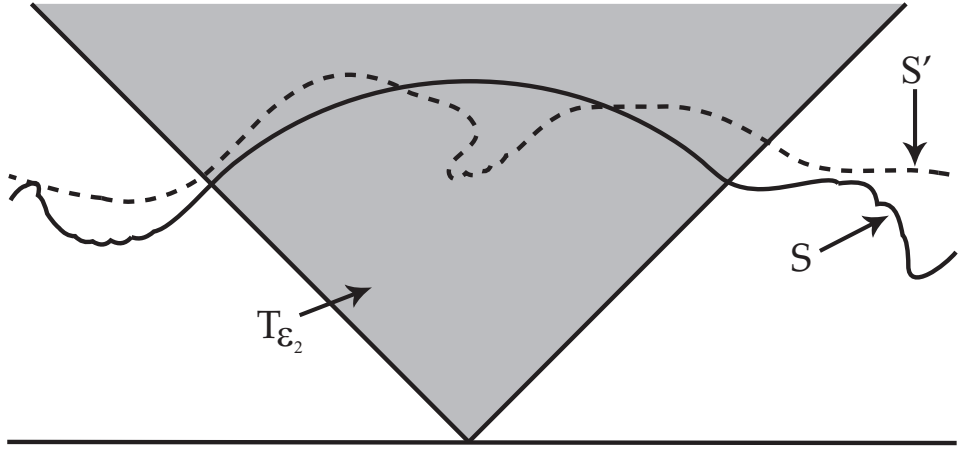}
\caption{ $S$ is nice inside $M_{(0,\epsilon_2]}$ and $S'$ is nice outside $M_{(0,\epsilon_2]}$.  We will join them by annuli, keeping the principal curvatures bounded.}
\end{figure}

Let $D$ be a totally geodesic disk in $S \cap M_{[\epsilon_2 ,\infty)}$.  We may assume that $D$ is transverse to the $2$-skeleton of $\mathcal{T}$.  We want to join the totally geodesic disk $D$ to the surface $S'$ consisting of geometric pieces by a surface with universally bounded principal curvatures.  We will replace some of the totally geodesic disks in $D$ by normal disks which connect $D$ to $S'$.  The idea is to pin down some of the vertices of a normal disk in $D$ and drag the others, along with some of the disk to meet the corresponding geometric piece.  In order to use a compactness argument to universally bound the principal curvatures of these ``connecting" disks, we may need to perturb $D$ first so that it does not come too close to vertices of $\mathcal{T}$ in the $1$-neighborhood $\mathcal{N}_1 (\bd M_{(0,\epsilon_1]})$ of $\bd M_{(0,\epsilon_1]}$.

We will perturb $D$ whenever it comes close to a vertex of $\mathcal{T}$ in $\mathcal{N} _1 (\bd M_{(0,\epsilon_1]})$ as follows.  Let $p \in \h^3$ and let $\mathcal{F}$ be a foliation of $\h^3$ by totally geodesic planes.  Replace the $(\mu / 100)$-ball around $p$ with the a smooth singular foliation $\mathcal{F}$ of the complement of a small region containing $p$ as shown in figure 9.  If $D$ intersects a $(\mu /100)$-ball centered at a vertex of $\mathcal{T}$, then we can consider it as a disk in a non-singular leaf of the foliation $\mathcal{F}$ (where the vertex is the point $p$), and replace it by a disk in a nonsingular leaf of $\mathcal{F'}$.  The principal curvatures of the leaves in $\mathcal{F'}$ are uniformly bounded.

\begin{figure}\label{fig9}
\includegraphics[width=\textwidth]{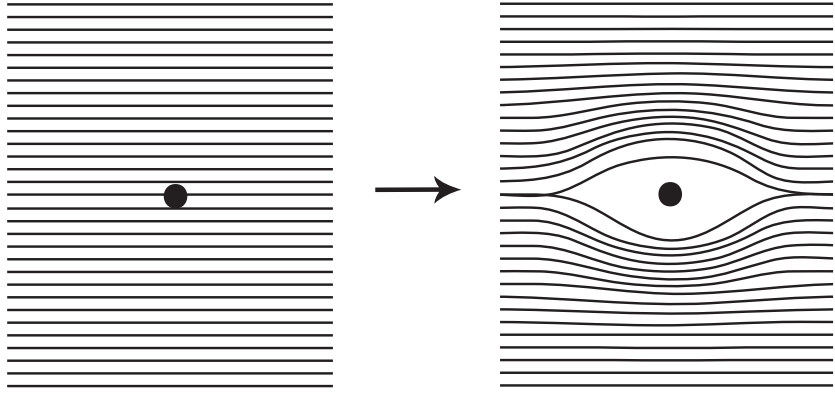}
\caption{ Cross-section view of how to change a foliation of $\h^3$ by totally geodesic planes near a vertex.}
\end{figure}

Let $t \in \mathcal{R} (L)$ and put $t$ inside $\h^3$.  Let $\mathcal{F}$ be a foliation of $\h^3$ by totally geodesic planes.  Change $\mathcal{F}$ near each vertex of $t$ as above to get a singular foliation $\mathcal{F'}$.  Let $d$ be the intersection of a non-singular leaf of $\mathcal{F'}$ with $t$.  For each proper subset of the set of vertices of $d$, choose a smooth normal disk which agrees with $d$ near these vertices and which agrees with the corresponding geometric disk near the other vertices.  Choose these disks so that they change continuously with the leaves of the foliation.  Call these new disks \textit{connecting disks}.  We can assure that the connecting disks meet to form a smooth surface when two tetrahedra are glued together by first defining them on a small neighborhood of each isometry class of triangles which appear as faces in $\mathcal{R} (L)$, and then extending them to the interiors of the tetrahedra.

The principal curvatures of the connecting disks in a tetrahedron $t$ corresponding a fixed foliation $\mathcal{F'}$ are uniformly bounded by compactness.  As the foliation $\mathcal{F'}$ in $\h^3$ changes via isotopy while $t$ remains fixed, isotop the connecting disks in $t$ continuously, so that compactness (of the set of foliations of $h^3$ by totally geodesic planes) implies that the set of possible connecting disks for $t$ has uniformly bounded principal curvatures.  Finally, the compactness of the set $\mathcal{R} (L)$ implies that the principal curvatures of any possible connecting disk in any tetrahedron in $\mathcal{R} (L)$ for any foliation $\mathcal{F'}$ are uniformly bounded in absolute value, say by $C(L)$.  Assume that $C(L) \ge B(L) + 1$.

Let $X$ be the union of the set of tetrahedra in $T$ which have at least one vertex in $M_{(0,\epsilon_1]}$ and at least one vertex outside $M_{(0,\epsilon_1]}$.  The union $X$ is homeomorphic to a torus times a compact interval.  Let $t$ be a tetrahedron in $X$ and replace each disk $d$ in $t \cap D$ as follows.  If $d$ has any vertices on the interior boundary of $X$ (i.e., the component of $\bd X$ which is closest to the center of $D$), then replace $d$ with the connecting disk which is totally geodesic near the vertices in the interior component of $X$ and which meets the corresponding geometric piece near the other vertices.  If $d$ has no vertices on the interior boundary of $X$, then replace $d$ with the corresponding geometric piece.  We end up with a disk $D'$ which agrees with $D$ near the interior boundary of $X$ and which agrees with $S'$ outside $M_{(0,\epsilon_1]}$, and which has principal curvatures bounded in absolute value by $C(L)$.

If $f: S \rightarrow M$ is an immersion, then the above argument works if we replace ``isotopy'' with ``homotopy''.\hfill$\Box$\\

\bibliographystyle{amsalpha}
\bibliography{tribound}

\newcommand{\etalchar}[1]{$^{#1}$}
\providecommand{\bysame}{\leavevmode\hbox to3em{\hrulefill}\thinspace}
\providecommand{\MR}{\relax\ifhmode\unskip\space\fi MR }
\providecommand{\MRhref}[2]{%
  \href{http://www.ams.org/mathscinet-getitem?mr=#1}{#2}
}
\providecommand{\href}[2]{#2}
\begin{thebibliography}{MTTW96}

\bibitem[Bre09]{Breslin2}
William Breslin, \emph{Thick triangulations of hyperbolic n-manifolds}, Pac. J.
  Math. \textbf{241} (2009), no.~2, 215--225.

\bibitem[Cas]{Casson}
A.~Casson, \emph{Three-dimensional topology}, Lecture notes from Beijing.

\bibitem[ELM{\etalchar{+}}00]{Edels}
Herbert Edelsbrunner, Xiang-Yang Li, Gary Miller, Andreas Stathopoulos, Dafna
  Talmor, Shang-Hua Teng, Alper Unger, and Noel Walkington, \emph{Smoothing and
  cleaning up slivers}, Proceedings of the Thirty-Second Annual ACM Symposium
  on Theory of Computing, ACM, 2000, pp.~273--277.

\bibitem[Fen89]{Fenchel}
Werner Fenchel, \emph{Elementary geometry in hyperbolic space}, de Gruyter
  Studies in Mathematics, no.~11, Walter de Gruyter \& Co., Berlin, 1989.

\bibitem[FHS83]{FHS}
Michael Freedman, Joel Hass, and Peter Scott, \emph{Least area incompressible
  surfaces in 3-manifolds}, Invent. Math. \textbf{71} (1983), 609--642.

\bibitem[Hak61]{Haken}
W.~Haken, \emph{Ein verfahren zur aufspaltung einer 3-mannigfaltigkeit in
  irreduzible 3 mannigfaltigkeiten}, Math. Zeit. \textbf{76} (1961), 427--467.

\bibitem[KM68]{Margulis}
D.~Kazhdan and G.~Margulis, \emph{A proof of {S}elberg's conjecture}, Math.
  USSR Sb. \textbf{4} (1968), 147--152.

\bibitem[Li00]{Li}
Xiang-Yang Li, \emph{Spacing control and sliver-free delaunay mesh},
  Proceedings, 9th International Meshing Roundtable, Sandia National
  Laboratories, 2000, pp.~295--306.

\bibitem[LL00]{Leibon}
Gregory Leibon and David Letscher, \emph{Delaunay triangulations andvoronoi
  diagrams for riemannian manifolds}, Proceedings of the Sixteenth Annual
  Symposium on Computational Geometry (Hong Kong, 2000) (2000), 341--349.

\bibitem[MTTW96]{Miller}
G.L. Miller, D.~Talmor, S.H. Teng, and H.~Walkington, N.and~Wang, \emph{Control
  volume meshes using sphere packing: generation refinement, and coarsening},
  Proceedings, 5th International Meshing Roundtable, Sandia National
  Laboratories, 1996, pp.~47--61.

\bibitem[Rub97]{Rub}
J.H. Rubinstein, \emph{Polyhedral minimal surfaces, {H}eegaard splittings and
  decision problems for 3-dimensional manifolds}, Geometric topology (Athens,
  GA, 1993), Amer. Math. Soc., 1997, pp.~1--20.

\bibitem[Sau05]{saucan3}
Emil Saucan, \emph{Note on a theorem of {M}unkres}, Mediterr. J. Math.
  \textbf{2} (2005), no.~2, 215--229.

\bibitem[Sau06a]{saucan1}
\bysame, \emph{The existence of quasimeromorphic mappings}, Ann. Acad. Sci.
  Fenn. Math. \textbf{31} (2006), no.~1, 131--142.

\bibitem[Sau06b]{saucan2}
\bysame, \emph{The existence of quasimeromorphic mappings in dimension 3},
  Conform. Geom. Dyn. \textbf{10} (2006), 21--40.

\bibitem[Sch83]{Schoen}
Richard Schoen, \emph{Estimates for stable minimal surfaces in
  three-dimensional manifolds}, Ann.of Math. Stud. \textbf{103} (1983),
  111--126.

\bibitem[Sto00]{Stocking}
M.~Stocking, \emph{Almost normal surfaces in 3-manifolds}, Trans. Amer. Math.
  Soc. \textbf{352} (2000), no.~1, 171--207.

\bibitem[SY79]{Schoen-Yau}
R.~Schoen and S.T. Yau, \emph{Existence of incompressible minimal surfaces and
  the topology of three dimensional manifolds with non-negative scalar
  curvature}, Annals of Math. \textbf{110} (1979), 127--142.

\end{thebibliography}

\end{document}